\documentclass[a4paper,11pt, reqno]{amsart}
\usepackage{amssymb,amsthm,amsmath}
\usepackage[T1]{fontenc}
\usepackage[colorlinks,citecolor=red,pagebackref,hypertexnames=false]{hyperref}

\numberwithin{equation}{section}
\allowdisplaybreaks[2]
\theoremstyle{plain}
\newtheorem{theorem}{Theorem}[section]

\newtheorem{lemma}[theorem]{Lemma}

\theoremstyle{definition}

\theoremstyle{remark}
\newtheorem{remark}[theorem]{Remark}

\newtheorem{case[theorem]}{Case}

\def\norm#1.#2.{\lVert#1\rVert_{#2}}

\author{Sunit Ghosh, Younis Ahmad Bhat and Jitendriya Swain}

	\address{Sunit Ghosh, Department of Mathematics, Indian Institute of Technology Guwahati, India.}
	\email{g.sunit@iitg.ac.in}
	
	\address{Younis Ahmad Bhat, Department of Mathematics, Indian Institute of Technology Guwahati, India.}
	\email{younis@iitg.ac.in}
	
	\address{Jitendriya Swain,   Department of Mathematics, Indian Institute of Technology Guwahati, India.}
	\email{jitumath@iitg.ac.in}

 \keywords{ Uncertainty principles, Dunkl operator, Dunkl transform, fractional Dunkl transform}

\subjclass[2010]{Primary  44A15, 42A38  ; Secondary 94A12.}

\date{\today}

\begin{document}
\baselineskip=15pt
\markboth{} {}

\bibliographystyle{amsplain}
\title[HPWUP for the fractional Dunkl transform on the real line]
{Heisenberg-Pauli-Weyl uncertainty principles for the fractional Dunkl transform on the real line }


\begin{abstract}
	 The aim of the paper is two-fold. First, we provide an explicit form of the functions for which equality holds for the uncertainty inequalities studied in \cite{Fei}. Second, we establish an $L^p$-type Heisenberg-Pauli-Weyl uncertainty principle for the fractional Dunkl transform, with $1 \leq p \leq 2$. For the case $p = 2$, we further derive a sharper uncertainty principle for the fractional Dunkl transform. Furthermore, we derive conditions leading to equality in both  the uncertainty principles obtained.
\end{abstract}

\maketitle
\def\BC{{\mathbb C}} \def\BQ{{\mathbb Q}}
\def\BR{{\mathbb R}} \def\BI{{\mathbb I}}
\def\BZ{{\mathbb Z}} \def\BD{{\mathbb D}}
\def\BP{{\mathbb P}} \def\BB{{\mathbb B}}
\def\BS{{\mathbb S}} \def\BH{{\mathbb H}}
\def\BE{{\mathbb E}}
\def\BN{{\mathbb N}}
\def\LP{{W(L^p(\BR^d, \BH), L^q_v)}}
\def\LPN{{W_{\BH}(L^p, L^q_v)}}
\def\LPQ{{W_{\BH}(L^{p'}, L^{q'}_{1/v})}}
\def\L1{{W_{\BH}(L^{\infty}, L^1_w)}}
\def\LB{{L^p(Q_{1/ \beta}, \BH)}}
\def\SP{S^{p,q}_{\tilde{v}}(\BH)}
\def\f{{\bf f}}
\def\h{{\bf h}}
\def\hp{{\bf h'}}
\def\m{{\bf m}}
\def\g{{\bf g}}
\def\ga{{\boldsymbol{\gamma}}}
\vspace{-.6cm}
\section{Introduction}
The classical Heisenberg-Pauli-Weyl uncertainty principle (HPWUP) states that for $f \in L^2(\mathbb{R})$ with $\|f\|_{L^2(\mathbb{R})} = 1$,  
\begin{align}\label{90}
	{\Delta}^2_{f,2} {\Delta}^2_{{\hat{f}},2}\geq\frac{1}{4},  
\end{align}
where ${\Delta}^2_{f,2}= \displaystyle\int_{-\infty}^{+\infty}|(x-( x)_{f})f(x)|^2dx, ( x)_{f}= \displaystyle\int_{-\infty}^{+\infty} x|f(x)|^2 dx $ and the classical Fourier transform $\hat{f}$ of $f$ defined by
\begin{align}
	\hat{f}(w)=\frac{1}{\sqrt{2\pi}}\int_{-\infty}^{+\infty} f(x) e^{-iwx} dx.
\end{align} 
In 1946, Gabor \cite{3} introduced the HPWUP in the context of signal analysis. Since then, numerous different forms of the HPWUP have emerged through various mathematical formulations \cite{CLF, FS, Havin}. It has significant applications in the field of signal processing, optics and so on(see \cite{dang,zayed, Zhao, 27, gxu}). Although inequality (\ref{90}) is the simplest form, it is not the most refined. Subsequently, Cohen \cite{cohen1, 1} and Dang-Deng-Qian \cite{dang} independently refined this classical result, deriving sharper uncertainty principles for complex-valued functions. An interesting extension was provided by Cowling-Price \cite{cowling}, who refined the classical HPWUP for $L^p$-type functions with $1\leq p \leq2$, yielding the following inequality
\begin{align}\label{9413}
	{\Delta}^2_{f,p} {\Delta}^2_{{\hat{f}},p}\geq\frac{1}{4},  
\end{align}
where ${\Delta}^p_{f,p}= \int_{-\infty}^{+\infty}|(x-\langle x\rangle_{f})f(x)|^pdx$. Later, Zhang \cite{Zhang} made further improvements to inequality (\ref{9413}).

The fractional Fourier transform is an extension of the classical Fourier Transform.  We refer to \cite{Namias, Kerr,zayed, Zhao} for its properties and applications. The inequalities (\ref{90}) and (\ref{9413}) are derived for the Fractional Fourier Transform for $L^2$-type HPWUP in \cite{Zhao} and for $L^p$-type HPWUP in \cite{Bahri} with $1\leq p\leq 2$. These inequalities were further improved in \cite{dang} and \cite{mm} by obtaining a sharper lower bound.

In this paper, we improve HPWUP results for the Dunkl transform \cite{Rosler, Fei} and the fractional Dunkl transform \cite{sami}, strengthening previously known bounds. The Dunkl transform and the Fractional Dunkl Transform generalize the Fourier transform and the Fractional Fourier Transform, respectively, by incorporating the symmetry of reflection groups associated with
the differential-difference operator, defined as \begin{align} \label{you}
T_{\mu}f(x)=f^{\prime}(x)+\left(\mu+\frac{1}{2}
\right)\frac{f(x)-f(-x)}{x},\,\, f\in\mathcal{C}^1(\mathbb{R}), \mu \geq -1/2,
\end{align} where $f^{\prime}$ denotes the classical derivative of $f$.
Such operators have been introduced by C. F. Dunkl in a series of papers \cite{Dunkl1, Dunkl2, Dunkl3}, which plays an essential role in Dunkl's generalization of spherical harmonics leading to further developments in harmonic analysis.

Let $\mu\geq-1/2$, for $1\leq p < \infty$, the space $L^p_{\mu}(\mathbb{R})$ consists of all complex functions on $\mathbb{R}$ such that $\|f\|_{\mu,p}=\left(\displaystyle\int\limits_{-\infty}^{+\infty} |f(x)|^p \left|x\right|^{2\mu+1} dx\right)^{\frac{1}{p}} < \infty$ and $L^{\infty}_{\mu}(\mathbb{R})$ is given in the usual way. For $f\in L^1_{\mu}(\mathbb{R})$, the Dunkl transform is given by
\begin{align}\label{1445}
	D_{\mu}f(w)=\frac{1}{2^{\mu+1}\Gamma(\mu+1)}\int\limits_{-\infty}^{+\infty}f(x)E_{\mu}(-iwx)\left|x\right|^{2\mu+1}dx,
\end{align}
where $E_{\mu}$ denotes the one dimensional Dunkl kernel, given by
\begin{align}
	E_{\mu}(z)=j_{\mu}(iz)+\frac{z}{2(\mu+1)}j_{\mu+1}(iz),
\end{align}
with the normalized spherical Bessel function $j_{\mu}(z)=\Gamma(\mu+1)\displaystyle \sum\limits_{n=0}^{\infty}\dfrac{(-1)^n(\frac{z}{2})^{2n}}{n!\Gamma{(n+\mu+1)}},\, z\in\mathbb{C}$.\\
Notice that when $\mu=-1/2$, we just have $E_{-1/2}(-iwx)=e^{-iwx}$, hence $D_{-1/2}$ coincides with the Fourier transform on $\mathbb{R}$.
For $f\in L^2_{\mu}(\mathbb{R}) \cap L^p_{\mu}(\mathbb{R})$, $1\leq p\leq 2$ we denote $\langle x \rangle_f=\displaystyle\int\limits_{-\infty}^{+\infty} x|f(x)|^2 \left|x\right|^{2\mu+1}dx$, ~and 
\begin{align}
\Delta_{\mu,p}(f)=\left(\int\limits_{-\infty}^{+\infty} \left|(x-\langle x \rangle_f)f(x)\right|^p\left|x\right|^{2\mu+1}dx\right)^{\frac{1}{p}}.
\end{align}
For a complex function $f$ defined on $\mathbb{R}$, we write its even and odd parts by $f_e(x)=\frac{f(x)+f(-x)}{2}$ and $f_o(x)=\frac{f(x)-f(-x)}{2}$ respectively.

 The HPWUP for the Dunkl transform, established by R\"{o}sler-Voit \cite{Rosler}, has attracted significant attention in the literature and is stated as follows: 
\begin{theorem}\cite{Rosler}.~\label{THM0}
Let $xf(x)\in L^2_{\mu}(\mathbb{R})~ \text{and}~~xD_{\mu}(f)(x)\in L^2_{\mu}(\mathbb{R})$ with $\|f\|_{\mu,2}=1$. Then 
\begin{align} \label{Rosler}
	\Delta_{\mu,2}^2(f)\Delta_{\mu,2}^2(D_{\mu}(f))\geq\left\{\left({\mu}+\frac{1}{2}\right)\left(\|f_e\|^2_{\mu,2}-\|f_0\|^2_{\mu,2}\right)+\frac{1}{2}\right\}^2, 
\end{align} 
equality holds in (\ref{Rosler}) if and only if  $f$ has the form $f(x)=de^{-\frac{1}{2\zeta}x^2} E_{\mu}(bx)$, with $b\in\mathbb{C}$ and $\zeta>0$.
\end{theorem}

Suppose that the complex function function $f$ is of the form $f(x) = \rho(x) e^{i \varphi(x)}$, where the classical derivatives $\rho^{\prime}(x)$ and $\varphi^{\prime}(x)$  exists for all $x\in \mathbb{R}$. The covariance and the absolute covariance of $f$ are defined and computed in \cite{Fei} as
\begin{align}
	Cov_{\mu}(f) = \int_{-\infty}^{+\infty} x \phi'(x) |x|^{2\mu +1} dx - \langle x \rangle_{f} \langle x \rangle_{D_{\mu}(f)},
\end{align}
and
\begin{align}\label{097}
	COV_{\mu}(f)&=\int_{-\infty}^{+\infty} \bigg| \left(x-\langle x \rangle_{f}\right) \biggl\{T_{\mu}(\varphi)(x) + \frac{\mu+\frac{1}{2}}{x}\bigg(\frac{\sin[\varphi(x)-\varphi(-x)] \rho(-x)}{\rho(x)}-(\varphi(x)-\varphi(-x))\bigg)\nonumber\\
	&\qquad\qquad-\langle x\rangle_{D_\mu(f)}\biggr\}\bigg| \rho^2(x) |x|^{2\mu +1} dx.
\end{align}
 In \cite{Fei}, Fei-Wang-Yang improved the inequality (\ref{Rosler}) by incorporating covariance and absolute covariance, leading to a strengthened formulation. The best lower bound they obtained is as follows.
\begin{theorem}\cite{Fei}. \label{fe}
	Assume that $f(x)=\rho(x)e^{i\varphi(x)}$ with $\|f(x)\|_{\mu,2}^2=1$ and the classical derivatives $\rho^{\prime}(x)$ and $\varphi^{\prime}(x)$  exists for all $x\in \mathbb{R}$, with $xf(x)\in L^2_{\mu}(\mathbb{R})~ \text{and}~~xD_{\mu}(f)(x)\in L^2_{\mu}(\mathbb{R})$. Then
	\begin{align}\label{110}
		\Delta_{\mu}^2(f)\Delta_{\mu}^2(D_{\mu}(f))\geq  \mathcal{A}^2 + COV_{\mu}^2(f),
	\end{align}	
where 
\begin{align}\label{567}
	\mathcal{A} = \left(\mu + \frac{1}{2}\right) \int_{-\infty}^{+\infty} \rho(x) \rho(-x) \cos [\varphi(x) - \varphi(-x)] |x|^{2\mu +1} dx + \frac{1}{2}.
\end{align}
	Moreover, if $\varphi^{\prime}(x)$ is continuous and $\rho$ is nonzero almost everywhere, then the equality in (\ref{110}) holds if and only if there exist $\zeta \in \mathbb{R}\setminus{\left\{0\right\}}$ and $\xi>0$ such that
	\begin{align}\label{117}
		(x-\langle x \rangle_{f})\rho(x)=-\zeta \bigg(T_{\mu}\rho(x)+\frac{\mu+\frac{1}{2}}{x}\rho(-x)\left(1-\cos[\varphi(x)-\varphi(-x)]\right)\bigg) 
	\end{align}
	and \begin{align}\label{118}
		\left|x-\langle x \rangle_{f}\right|=\xi \bigg|T_{\mu}\varphi(x)+\frac{\mu+\frac{1}{2}}{x}\left(\frac{\sin[\varphi(x)-\varphi(-x)]}{\rho(x)}\rho(-x)-(\varphi(x)-\varphi(-x))\right)-\langle x\rangle_{D_\mu(f)}\bigg|.  
	\end{align}
\end{theorem}
\begin{remark}
	(1) If $f$ is of the form $f(x) = \rho(x)e^{i \varphi(x)}$, then $\mathcal{A}^2$ coincides with the right-hand side of (\ref{Rosler}). (See Section \ref{sec2}). \\
	 (2) For the equality in (\ref{110}), the authors overlooked the fact that $(x-\langle x \rangle_{f})\rho(x)$ and $T_{\mu}\rho(x)+\frac{\mu+\frac{1}{2}}{x}\rho(-x)\left(1-\cos[\varphi(x)-\varphi(-x)]\right)$ must have the same or opposite sign. This follows from the second-to-last line on page 891 of \cite{Fei}. Combining this with (3.18) of \cite{Fei} yields (\ref{117}). 
\end{remark}

Let $\alpha \in \mathbb{R}$. For $f \in L^1_\mu(\mathbb{R})$, we consider the fractional Dunkl transform $D_\mu^\alpha f$ defined as
\begin{align}\label{2001}
D_{\mu}^{\alpha}f(\omega)=\left\{\begin{array}{ccc}N_{\mu,n}\int\limits_{-\infty}^{+\infty}  e^{\frac{i}{2}(x^2+{\omega}^2)\cot\alpha}E_{\mu}\Big(\frac{-\,i{\omega}x}{\sin \alpha}\Big)f(x) |x|^{2\mu+1} dx,& (2n-1)\pi<x<(2n+1)\pi,\\
	f(x),\qquad\qquad\qquad\qquad\qquad\qquad\qquad\qquad\qquad& x=2n\pi,\qquad\qquad\qquad\qquad\\
	f(-x),\qquad\qquad\qquad\qquad\qquad\qquad\qquad\qquad\qquad& x=(2n+1)\pi,\qquad\qquad\qquad
\end{array}\right. 
\end{align}
where $$N_{\mu,n}=\dfrac{e^{i(\mu+1)(\frac{\hat{\alpha}\pi}{2}-(\alpha-2n\pi))}}{\Gamma(\mu+1)(2|\sin\alpha|)^{\mu+1}}, ~~~~\mbox{with} ~~~~\hat{\alpha}=sgn(sin(\alpha)).$$
The fractional Dunkl transform provides a natural generalization of the Dunkl transform, see Section \ref{sec2} for more details.
In \cite{sami}, Ghazouani and Bouzeffour established the following uncertainty principle for the fractional Dunkl transform.
\begin{theorem} \cite{sami}.\label{01}
	Suppose $\alpha, \beta \in\mathbb{R}$ such that $\beta -\alpha \notin \pi\mathbb{Z}$, and let $x D_{\mu}^{\alpha}f(x)\in L^2_{\mu}(\mathbb{R})$ and  $x D_{\mu}^{\beta}f(x)\in L^2_{\mu}(\mathbb{R})$ with $\|f\|_{\mu,2}=1$. Then
	\begin{align}\label{sam}
		\Delta_{\mu, 2}^2( D_{\mu}^{\alpha}f)\Delta_{\mu,2}^2(D_{\mu}^{\beta}(f))&\geq \sin^2(\alpha-\beta)\left\{\left({\mu}+\frac{1}{2}\right)\left(\|f_e\|^2_{\mu,2}-\|f_0\|^2_{\mu,2}\right)+\frac{1}{2}\right\}^2.
	\end{align}
	Moreover, equality holds if and only if $f(x)=d e^{-\frac{1}{2\zeta}x^2}E_{\mu}(bx)$ with $b,d\in\mathbb{C}$ and $\zeta>0$.
\end{theorem}

Our main results are stated as follows. We solve the differential equations (\ref{117}) and (\ref{118}) and, as a result, provide an explicit form of the functions that attain equality in (\ref{110}), thereby completing the results in Theorem 3.3 of \cite{Fei}.

\begin{theorem} \label{THM1}
Assume that $f(x)=\rho(x)e^{i\varphi(x)}$ with $\|f(x)\|_{\mu,2}^2=1$ and the classical derivatives $\rho^{\prime}(x)$ and $\varphi^{\prime}(x)$  exists for all $ x\in\mathbb{R}$, with $\varphi^{\prime}(x)$ continuous and $\rho$ nonzero almost everywhere.  Then (\ref{117}) and (\ref{118}) holds if and only if $f$ has one of the following form
\begin{align}\label{12}
	f(x)=d_1 e^{-\frac{1}{2\zeta}x^2} e^{\frac{i}{2\xi}x^2} 	E_{\mu}(bx), 
\end{align}
or\begin{align}\label{13}
	f(x)=d_2 e^{-\frac{1}{2\zeta}x^2} e^{-\frac{i}{2\xi}x^2} 	E_{\mu}(b^{\prime}x), 
\end{align}
or\begin{align}\label{14}
	f(x)=\left\{\begin{array}{cc} d_3 e^{-\frac{1}{2\zeta}x^2} e^{\frac{i}{2\xi}x^2} 	E_{\mu}(bx),& x\geq \langle x\rangle_{f}, \\
		\quad d_4 e^{-\frac{1}{2\zeta}x^2} e^{-\frac{i}{2\xi}x^2} 	E_{\mu}(b^{\prime}x),& x <\langle x\rangle_{f}, \\
	\end{array}\right. 
\end{align}
or
\begin{align} \label{15}
	f(x)=\left\{\begin{array}{cc} \quad d_5 e^{-\frac{1}{2\zeta}x^2} e^{-\frac{i}{2\xi}x^2} 	E_{\mu}(b^{\prime}x),& x\geq \langle x\rangle_{f}, \\
		d_6 e^{-\frac{1}{2\zeta}x^2} e^{ \frac{i}{2\xi}x^2} 	E_{\mu}(bx),& x <\langle x\rangle_{f}, \\
	\end{array}\right. 
\end{align}
where $b=\left(\frac{1}{\zeta}-\frac{i}{\xi}\right)\langle x\rangle_{f}+i \langle x\rangle_{D_{\mu}(f)}$ and $b^{\prime}=\left(\frac{1}{\zeta}+\frac{i}{\xi}\right)\langle x\rangle_{f} +i \langle x\rangle_{D_{\mu}(f)}$ for some $\zeta, \xi >0$ and $d_1, d_2, d_3, d_4, d_5, d_6\in\mathbb{C}$ chosen such that $\|f\|_{\mu,2}=1$.
\end{theorem}  
For the fractional Dunkl transform, we obtain the following $L^p-$type HPWUP, $1 \leq p \leq 2$, in more general situation. 
\begin{theorem} \label{MTHM1}
	Assume that $f(x)=\rho(x)e^{i\varphi(x)}$ with $\|f(x)\|_{\mu,2}^2=1$ and the derivatives $\rho^{\prime}(x)$, $\varphi^{\prime}(x)$ and $f^{\prime}(x)$ exists for all $ x\in\mathbb{R}$. Let $\alpha, \beta \in \mathbb{R}$ be such that  $\beta - \alpha \notin \pi\mathbb{Z}$ and let $x D_{\mu}^{\alpha}f(x)\in L^2_{\mu}(\mathbb{R})\cap L^p_{\mu}(\mathbb{R})$, $x D_{\mu}^{\beta}f(x)\in L^p_{\mu}(\mathbb{R}), 1\leq p\leq 2 $, then 
	\begin{align} \label{16}
		\notag\Delta_{\mu, p}^2( D_{\mu}^{\alpha}f)\Delta_{\mu,p}^2(D_{\mu}^{\beta}(f&))\geq \frac{|\sin(\beta-\alpha)|^{2(\mu +1)(\frac{2}{p}-1)}}{\left(2^{\mu+1}\Gamma(\mu +1)\right)^{2(\frac{2}{p}-1)}}\Biggl\{\sin^2(\beta-\alpha)\left((\mu+\frac{1}{2})(\|f_e\|^2_{\mu,2}-\|f_o\|^2_{\mu,2})+\frac{1}{2}\right)^2\\ 
		&+\bigg(\cos\alpha\,\cos\beta \,\Delta_{\mu, 2} ^2 (f)+\sin(\alpha+\beta)\, Cov_{\mu}(f)+\sin\alpha \sin\beta\, \Delta_{\mu,2}^2(D_{\mu}(f))\bigg)^2\Biggr\}. 
	\end{align}
	Moreover, if $\varphi^{\prime}(x)$ is continuous and $\rho$ is nonzero almost everywhere, then the equality in (\ref{16}) holds if and only if $p = 2$ and  $f$ is of the form (\ref{12}) or (\ref{13}).
\end{theorem}
\begin{remark} (1) For $p=2$, Theorem \ref{MTHM1} provides a stronger lower bound than the one presented in equation (\ref{sam}) (see Theorem 6.1 of \cite{sami}).\\
(2)  When $\alpha =0, \beta =\pi/2$ and  $p=2$, Theorem \ref{MTHM1} coincides with Theorem 3.2 of \cite{Fei} and provides a necessary and sufficient condition for the uncertainty inequality to hold with equality.
\end{remark}

For $p=2$, we further improve the uncertainty inequality in (\ref{16}) by deriving a truly sharp lower bound, which is  established in the following theorem.
\begin{theorem} \label{MTHM2}
	Assume that $f(x)=\rho(x)e^{i\varphi(x)}$ with $\|f(x)\|_{\mu,2}^2=1$. Let $\alpha, \beta \in \mathbb{R}$ be such that  $\beta - \alpha \notin \pi\mathbb{Z}$, and let $x D_{\mu}^{\alpha}f(x)\in L^2_{\mu}(\mathbb{R})$ and  $x D_{\mu}^{\beta}f(x)\in L^2_{\mu}(\mathbb{R})$, then
	\begin{align} \label{17}
		\notag \Delta_{\mu, 2}^2( D_{\mu}^{\alpha}f)\Delta_{\mu,2}^2(D_{\mu}^{\beta}(f))&\geq \sin^2(\alpha-\beta)\Biggl(\left\{\left({\mu}+\frac{1}{2}\right)\left(\|f_e\|^2_{\mu,2}-\|f_0\|^2_{\mu,2}\right)+\frac{1}{2}\right\}^2+COV_{\mu}^2(f)\\
		-Cov_{\mu}(f)\Biggl)&+\biggl(\cos\alpha \, \cos\beta \Delta_{\mu, 2}^2 (f) +\sin(\alpha+\beta)Cov_{\mu}(f)+\sin\alpha\,\sin\beta \Delta_{\mu,2}^2(D_{\mu}(f))\biggr)^2.
	\end{align}
	Moreover, if $\varphi^{\prime}(x)$ is continuous and $\rho$ is nonzero almost everywhere, then the equality in (\ref{17}) holds if and only if $f$ is of the form (\ref{12}), (\ref{13}), (\ref{14}) or (\ref{15}).
\end{theorem}
\begin{remark}
	For $\mu = -\frac{1}{2}$, the results presented in Theorems \ref{MTHM1} and \ref{MTHM2} are obtained by Chen-Dang-Mai \cite{Chen} and Dang-Deng-Qian \cite{mm}, respectively.  
\end{remark}

The paper is structured as follows:  In Section \ref{sec2}, we provide some essential results related to the Dunkl transform and the fractional Dunkl transform. In Section \ref{sec3}, we establish our main results, followed by a concluding remark. 

\section{Preliminaries}\label{sec2}
In this section, we outline some important results related to the Dunkl transform and the fractional Dunkl transform.
\subsection{Dunkl transform and the Dunkl operator}
For any function $f\in L^1_{\mu}(\mathbb{R})$, the Dunkl transform of $f$ defined in (\ref{1445}), $D_\mu(f)$ belongs to $  C_0(\mathbb{R})$ and satisfies
\begin{align}
	\|D_\mu(f)\|_{\mu,\infty} \leq \frac{1}{2^{\mu+1}\Gamma(\mu +1)}\|f\|_{\mu,1}.
\end{align}
If $f \in L^2_\mu(\mathbb{R})$, then $D_\mu(f)$ belongs to $ L^2_\mu(\mathbb{R})$ and satisfies the following Plancherel formula
\begin{align}
	\|D_{\mu}(f)\|_{\mu,2}=\|f\|_{\mu,2}.
\end{align}
Interpolating these two inequalities gives the Hausdorff-Young inequality: if \( f \in L^p_{\mu}(\mathbb{R}) \), where \( 1 \leq p \leq 2 \) and \( \frac{1}{p} + \frac{1}{q} = 1 \), then \( D_\mu(f) \in L^q_{\mu}(\mathbb{R}) \), and 
\begin{align}\label{HYin}
	\|D_\mu(f)\|_q \leq \frac{1}{(2^{\mu+1}\Gamma(\mu +1))^{\frac{2}{p}-1}}\|f\|_p,
\end{align}
with equality if and only if \( p = q = 2 \).

\noindent Suppose $f, g \in L^2_\mu(\mathbb{R})$ such that $T_\mu f, T_\mu g \in L^2_\mu(\mathbb{R})$, then we have
\begin{align}
	\int_{-\infty}^{+\infty} T_\mu f(x) g(x) |x|^{2\mu +1} dx = -\int_{-\infty}^{+\infty}  f(x) T_\mu g(x) |x|^{2\mu +1}dx
\end{align}
and 
\begin{align}
	D_\mu (T_\mu f)(w) = iw D_\mu(f)(w).
\end{align}
For $f,g \in C^1(\mathbb{R}\setminus\{0\})$ the following product formula holds
\begin{align}
	T_\mu(fg) = (T_\mu f)\cdot g + f\cdot (T_\mu g).
\end{align}
We refer to \cite{Dunkl3, jeu, Rosler1, Rosler} for a detailed study on Dunkl transform and its properties.

\subsection{Fractional Dunkl transform} 
The fractional Dunkl transform introduced by Ghazouani-Bouzeffour \cite{sami} and defined in (\ref{2001}), was initially defined on \( L^1_\mu(\mathbb{R}) \) and admits a unique extension as a unitary operator on \( L^2_\mu(\mathbb{R}) \). If we denote this extension by \( D_\mu^{\alpha} \), the family \( \{D_\mu^\alpha\}_{\alpha \in \mathbb{R}} \) forms a group structure. Specifically, for any \( f \in L^2_\mu(\mathbb{R}) \), we have the following property:
\begin{align}
	D_{\mu}^{\alpha} \circ D_{\mu}^{\beta}(f)=D_{\mu}^{\alpha+\beta}(f)
\end{align}
with \( D^0_\mu \) is the identity, and the inverse is given by \( (D_\mu^\alpha)^{-1} = D_\mu^{-\alpha} \). The fractional Dunkl transform \(D_\mu^{\alpha} \) also has the following exponential form. 
$$D_\mu^{\alpha} = e^{-i\alpha(\mu+1)}e^{-i\frac{\alpha}{2}(T_\mu^2-x^2)}.$$

Let $\alpha \notin \pi \mathbb{Z}$. The following relation holds between the fractional Dunkl transform and the Dunkl transform 
 \begin{align} \label{36}
	D^{\alpha}_{\mu}(f)(w)=\frac{e^{i(\mu+1)(\frac{\hat{\alpha}\pi}{2}-(\alpha-2n\pi))}}{\left|\sin\alpha\right|^{\mu+1}}e^{\frac{iw^2\cot\alpha}{2}}D_{\mu}(g)\left(\frac{w}{\sin\alpha}\right).
\end{align}
where $g(x) = e^{\frac{ix^2\cot\alpha}{2}} f(x)$ and $\hat{\alpha}=sgn(sin(\alpha))$.

\subsection{Covariance and absolute covariance}\label{sub2.3} The covariance of $f(x)=\rho(x)e^{i\varphi(x)}$ can also be expressed as follows (see \cite{Fei}):
\begin{align}\label{098}
	Cov_{\mu}(f)&=\int_{-\infty}^{+\infty} \left(x-\langle x \rangle_{f}\right) \biggl\{T_{\mu}(\varphi)(x) + \frac{\mu+\frac{1}{2}}{x}\bigg(\frac{\sin[\varphi(x)-\varphi(-x)] \rho(-x)}{\rho(x)}-(\varphi(x)-\varphi(-x))\bigg)\nonumber\\
	&\qquad\qquad-\langle x\rangle_{D_\mu(f)}\biggr\} \rho^2(x) |x|^{2\mu +1} dx.
\end{align}
From the fact  $\displaystyle\int_{-\infty}^{+\infty} |f(x)| |x|^{2\mu+1} dx \geq \displaystyle\int_{-\infty}^{+\infty} f(x) |x|^{2\mu+1} dx $, it follows that $COV_\mu(f) \geq Cov_\mu(f)$. Assuming that $\varphi'(x)$ is continuous, it clear from equations (\ref{098}) and (\ref{097}) that $COV_\mu(f) = Cov_\mu(f)$ if and only if $(x-\langle x \rangle_{f})$ and the expression $$T_{\mu}(\varphi)(x) + \frac{\mu+\frac{1}{2}}{x}\bigg(\frac{\sin[\varphi(x)-\varphi(-x)] \rho(-x)}{\rho(x)}-(\varphi(x)-\varphi(-x))\bigg)
-\langle x\rangle_{D_\mu(f)}$$ have the same or opposite signs.

Suppose  $xf(x)\in L^2_{\mu}(\mathbb{R})~ \text{and}~~xD_{\mu}(f)(x)\in L^2_{\mu}(\mathbb{R})$ with $\|f\|_{\mu,2}=1$. Then we have the following relation (see \cite{Rosler}):
\begin{align}\label{yt4}
Re \int_{-\infty}^{+\infty} x f(x)\overline{T_\mu f(x)} |x|^{2\mu+1} dx =	-\left(\|f_e\|^2_{\mu,2}-\|f_0\|^2_{\mu,2}\right)-\frac{1}{2},
\end{align} where left hand side denotes the real part of the integral.
Also we can write
\begin{align}\label{tr4}
	\int_{-\infty}^{+\infty} x f(x)\overline{T_\mu f(x)} |x|^{2\mu+1} dx = \frac{1}{2} \int_{-\infty}^{+\infty}  &f(x)\overline{(xT_\mu - T_\mu x) f(x)} |x|^{2\mu+1} dx\nonumber\\ &+ \frac{1}{2} \int_{-\infty}^{+\infty}  f(x)\overline{(xT_\mu + T_\mu x) f(x)} |x|^{2\mu+1} dx.
\end{align}
A straightforward calculation (see the proof of Theorem 3.2 in \cite{Fei}) yields,
$$(xT_\mu - T_\mu x) f(x) = - f(x) - (2\mu +1) f(-x),$$ and $$(xT_\mu + T_\mu x) f(x) = 2 xf'(x) + 2(\mu +1)f(x). $$
Now, assume that $f(x) = \rho(x) e^{i \varphi(x)}$, where the derivatives $\rho^{\prime}(x)$ and $\varphi^{\prime}(x)$  exists for all $x\in \mathbb{R}$, then, from (\ref{tr4}), we obtain 
\begin{align}\label{re4}
\int_{-\infty}^{+\infty} x f(x)\overline{T_\mu f(x)} |x|^{2\mu+1} dx = -\mathcal{A} - i \left(Cov_\mu(f) + \langle x\rangle_{f} \langle x\rangle_{D_\mu(f)}\right),
\end{align}
where $\mathcal{A}$ is defined in (\ref{567}). By comparing (\ref{yt4}) and (\ref{re4}), we conclude that
\begin{align}\label{uyt}
	\mathcal{A} = \left(\|f_e\|^2_{\mu,2}-\|f_0\|^2_{\mu,2}\right) + \frac{1}{2}.
\end{align}

\section{HPWUP for the fractional Dunkl Transform:}\label{sec3}
Before proving our main results, we start with the following lemmas.
 \begin{lemma}\cite{Fei} \label{lem1} Let $f(x)=\rho(x)e^{i\varphi(x)}\in L^2_{\mu}(\mathbb{R})$ with $\|f(x)\|_{\mu,2}=1$. Suppose that the classical derivatives $\rho^{\prime}(x)$, $\varphi^{\prime}(x)$ exists for all $x\in\mathbb{R}$ and $f^{\prime}(x)\in L^2_{\mu}(\mathbb{R})$, then 
{\small
	\begin{align}
	\notag	&\Delta_{\mu,2}^2(D_{\mu}(g))=\int\limits_{-\infty}^{+\infty} \left\{T_{\mu}\rho(x)+\frac{\mu+\frac{1}{2}}{x}\rho(-x)\big(1-\cos[\varphi(x)-\varphi(-x)]\big)\right\}^2|x|^{2\mu+1} dx\\
		&+ \int\limits_{-\infty}^{+\infty}\left\{T_{\mu}\varphi(x)+\frac{\mu+\frac{1}{2}}{x}\biggl(\frac{\sin[\varphi(x)-\varphi(-x)]}{\rho(x)}\rho(-x)-(\varphi(x)-\varphi(-x))\biggr)-\langle x\rangle_{D_\mu(f)}\right\}^2\rho^2(x)|x|^{2\mu+1} dx.
\end{align}}
 \end{lemma}
\begin{lemma}
 \label{lem2} Let $f(x)=\rho(x)e^{i\varphi(x)}$ with $\|f(x)\|_{\mu,2}=1$ and $g(x)=e^{i\frac{x^2\cot\alpha}{2}}f(x); \, \alpha\in\mathbb{R}\setminus\pi\mathbb{Z}$. Suppose that the classical derivatives $\rho^{\prime}(x)$, $\varphi^{\prime}(x)$ exists for all $ x\in\mathbb{R}$ and $f^{\prime}(x)\in L^2_{\mu}(\mathbb{R})$, then 
 \begin{align}
	\langle x\rangle_{D_{\mu}(g)} = \cot\alpha\,\langle x\rangle_f+\langle x\rangle_{D_{\mu}(f)},
\end{align}
and 
\begin{align}
	\Delta_{\mu,2}^2(D_{\mu}(g)) = \Delta_{\mu,2}^2(D_{\mu}(f)) + 2\cot\alpha ~ Cov_{\mu}(f) + \cot\alpha\,\Delta_{\mu,2}^2(f).
\end{align}
\end{lemma}
\noindent{\bf{Proof}}.	We have
	\begin{align*}
		\langle x\rangle_{D_{\mu}(g)}=\int\limits_{-\infty}^{+\infty} x D_{\mu}(g)(x) \overline{D_{\mu}(g)(x)} |x|^{2\mu+1} dx=\int\limits_{-\infty}^{+\infty} -i T_{\mu}(g)(x)\overline{g(x)} |x|^{2\mu+1} dx.
	\end{align*}
	On substituting $T_{\mu}(g)(x)=e^{i\frac{x^2\cot\alpha}{2}}\left(ix\cot\alpha\, f(x)+T_{\mu}f(x)\right)$, we get
	\begin{align*}
		\langle x\rangle_{D_{\mu}(g)}&=\cot\alpha\, \int\limits_{-\infty}^{+\infty} x\left|f(x)\right|^2 \left|x\right|^{2\mu+1} dx +\int\limits_{-\infty}^{+\infty} -i T_{\mu}f(x)\overline{f(x)}  \left|x\right|^{2\mu+1} dx\\
		&=\cot\alpha \, \langle x\rangle_f+\langle x\rangle_{D_{\mu}(f)}.
	\end{align*}
	Applying Lemma \ref{lem1}, we get
	\begin{align*}
		\Delta_{\mu,2}^2(D_{\mu}(g))&=\int\limits_{-\infty}^{+\infty} \left\{T_{\mu}\rho(x)+\frac{\mu+\frac{1}{2}}{x}\rho(-x)\left(1-\cos[\varphi(x)-\varphi(-x)]\right)\right\}^2 \rho^2(x) \left|x\right|^{2\mu+1}dx\\
		&+\int\limits_{-\infty}^{+\infty} \biggl\{\varphi^{\prime}(x)+x\cot\alpha+\frac{\mu+\frac{1}{2}}{x}\left(\frac{\sin[\varphi(x)-\varphi(-x)]}{\rho(x)}\rho(-x)\right)-\cot\alpha \, \langle x \rangle_{f} \\
		&\qquad\qquad -\langle x \rangle_{D_{\mu}(f)}\biggr\}^2 \rho^2(x) \left|x\right|^{2\mu+1}dx\\
		&=\Delta_{\mu,2}^2(D_{\mu}(f))+\cot^2\alpha \, \Delta_{\mu,2}^2(f)+2\cot\alpha \, \int\limits_{-\infty}^{+\infty}\left(x- \langle x \rangle_{f}\right)\\
		&\qquad\times\left\{\varphi^{\prime}(x)+\frac{\mu+\frac{1}{2}}{x}\left(\frac{\sin[\varphi(x)-\varphi(-x)]}{\rho(x)}\rho(-x)\right)-\langle x \rangle_{D_{\mu}(f)}\right\}\rho^2(x) \left|x\right|^{2\mu+1}dx\\
		&=\Delta_{\mu,2}^2(D_{\mu}(f))+\cot\alpha \,\Delta_{\mu,2}^2(f)+2\cot\alpha \,\left(\int\limits_{-\infty}^{+\infty} x\varphi^{\prime}(x)\rho^2(x)|x|^{2\mu+1}dx-\langle x \rangle_{f}\langle x \rangle_{D_{\mu}(f)}\right)\\
		&=\Delta_{\mu,2}^2(D_{\mu}(f))+2\cot\alpha~Cov_{\mu}(f)+\cot\alpha \,\Delta_{\mu,2}^2(f).
	\end{align*}
$\hfill\square$
\begin{lemma} \label{lem3}
 ~Let $\alpha\in\mathbb{R}$ and let $f(x)=\rho(x)e^{i\varphi(x)}$ with $\|f(x)\|_{\mu,2}=1$. Suppose that the classical derivatives $\rho^{\prime}(x)$, $\varphi^{\prime}(x)$ and $f^{\prime}(x)$ exists for all $ x\in\mathbb{R}$. If $xf(x)\in L^2_{\mu}(\mathbb{R})$ and $xD_{\mu}^{\alpha}f(x) \in L^2_{\mu}(\mathbb{R})$, then 
\begin{align}
	\langle x\rangle_{D_{\mu}^{\alpha}(f)}=\cos\alpha \, \langle x\rangle_{f}+\sin\alpha \, \langle x\rangle_{D_{\mu}(f)} 
\end{align}
and
\begin{align}
	\Delta_{\mu,2}^2(D_{\mu}^{\alpha}(f))&=\cos^2\alpha \, \Delta_{\mu,2}^2(f)+2\cos\alpha\, \sin\alpha\,  Cov_{\mu}(f)+\sin^2\alpha\, \Delta_{\mu,2}^2(D_{\mu}(f)).
\end{align}
\end{lemma}
\noindent{\bf{Proof}}. ~Using ({\ref{36}}), we can write
	\begin{align} \label{37}
	\langle x\rangle_{D_{\mu}^{\alpha}(f)}=\frac{1}{|\sin\alpha|^{(\mu+1)2}}\int\limits_{-\infty}^{+\infty} x\left|D_{\mu}(g)\left(\frac{x}{\sin\alpha}\right)\right| \left|x\right|^{2\mu+1}dx=\sin\alpha \, \langle x\rangle_{D_{\mu}(g)}, 
	\end{align}
   Where $g(x)= e^{\frac{ix^2 \cot \alpha}{2}} f(x)$. Using Lemma \ref{lem2}, we get
	\begin{align*}
		\langle x\rangle_{D_{\mu}^{\alpha}(f)}=\cos\alpha \, \langle x\rangle_{f}+\sin\alpha \, \langle x\rangle_{D_{\mu}(f)}.
	\end{align*}
	By using (\ref{36}) and (\ref{37}), we have
	\begin{align*}
		\Delta_{\mu,2}^2(D_{\mu}^{\alpha}(f))&=\int\limits_{-\infty}^{+\infty}\left(x-\langle x\rangle_{D_{\mu}^{\alpha}(f)}\right)^2 \left|D_{\mu}^{\alpha}(f)(x)\right|^2 \left|x\right|^{2\mu+1}dx\\
		&=\int\limits_{-\infty}^{+\infty} \left(x \sin\alpha-\langle x\rangle_{D_{\mu}^{\alpha}(f)}\right)^2 \left|D_{\mu}(g)(x)\right|^2 \left|x\right|^{2\mu+1}dx\\
		&=\sin^2\alpha \,\,  	\Delta_{\mu,2}^2(D_{\mu}(g)).
	\end{align*}
	Therefore, Lemma \ref{lem3} follows from Lemma \ref{lem2}.
$\hfill\square$\\

Now we proceed to prove Theorem \ref{THM1}. \vspace{.2cm}

\noindent{\bf{Proof of Theorem \ref{THM1}}}. Let $A(x)=T_{\mu}\varphi(x)+\frac{\mu+\frac{1}{2}}{x}\left(\frac{\sin[\varphi(x)-\varphi(-x)]}{\rho(x)}\rho(-x)-(\varphi(x)-\varphi(-x))\right)\, \forall\, x\in\mathbb{R}$. As the left-hand side of equation (\ref{118}) corresponds to the absolute value of a linear function, there are four possible scenarios:
\begin{align}\label{38}
 A(x)=\frac{1}{\xi} (x-\langle x\rangle_{f})+\langle x\rangle_{D_\mu(f)},
\end{align}
or
\begin{align}\label{39}
 A(x)=-\frac{1}{\xi} (x-\langle x\rangle_{f})+\langle x\rangle_{D_\mu(f)},
\end{align}
or
\begin{align}\label{40}
	A(x)=\left\{\begin{array}{cc}  \frac{1}{\xi} (x-\langle x\rangle_{f})+\langle x\rangle_{D_\mu(f)},& x\geq \langle x\rangle_{f}, \\
		-\frac{1}{\xi} (x-\langle x\rangle_{f})+\langle x\rangle_{D_\mu(f)},& x <\langle x\rangle_{f},\\
	\end{array}\right.
\end{align}
or
\begin{align}\label{41}
	A(x)=\left\{\begin{array}{cc}  -\frac{1}{\xi} (x-\langle x\rangle_{\mu})+\langle x\rangle_{D_\mu(f)},& x\geq \langle x\rangle_{f}, \\
		\frac{1}{\xi} (x-\langle x\rangle_{\mu})+\langle x\rangle_{D_\mu(f)},& x <\langle x\rangle_{f}.
	\end{array}\right.
\end{align}
First we assume that 
\begin{align}\label{ss}
A(x)=\frac{1}{\xi} (x-\langle x\rangle_{\mu})+\langle x\rangle_{D_\mu(f)}.
\end{align}
It follows that $f\in\mathcal{C}^{\infty}(\mathbb{R}\setminus{\left\{0\right\}})$ solves (\ref{117}) and (\ref{ss}), if $r(x) e^ {i\phi(x)} = F(x)=e^{\frac{1}{2\zeta}x^2}e^{\frac{i}{2\xi}x^2} f(x)$ satisfies
\begin{align}
 T_{\mu}r(x)+\frac{\mu+\frac{1}{2}}{x}r(-x)\left(1-\cos[\phi(x)-\phi(-x)]\right)=\frac{1}{\zeta} \langle x\rangle_{f} \,r(x)
\end{align}
and
\begin{align}
	 T_{\mu}\phi(x)+\frac{\mu+\frac{1}{2}}{x}\left(\frac{r(-x) \sin[\phi(x)-\phi(-x)]\,}{r(x)}-(\phi(x)-\phi(-x))\right)=-\frac{\langle x\rangle_{f}}{\xi} +\langle x\rangle_{D_\mu(f)},
\end{align}
which is of the form 
\begin{align}\label{314}
T_{\mu}F=bF, \, \, \, \, \text{where}\quad b=\left(\frac{1}{\zeta}-\frac{i}{\xi}\right)\langle x\rangle_{f}+i \langle x\rangle_{D_{\mu}(f)}.
\end{align}
In \cite{Rosler}, the authors derived the explicit solution to the equation (\ref{314}) for $F\in\mathcal{C}^{\infty}(\mathbb{R}\setminus{\left\{0\right\}})$ which is given by
\begin{align*}
F(x)=d_1 E_{\mu}(bx), \, \, \text{for some}\, d_1\in \mathbb{C}.
\end{align*}
Therefore $f(x)=d_1 e^{-\frac{1}{2\zeta}x^2}e^{\frac{i}{2\xi}x^2}E_{\mu}(bx)$. Now, if $\zeta<0$, then $f(x)=d_1 e^{-\frac{1}{2\zeta}x^2}e^{\frac{i}{2\xi}x^2}E_{\mu}(bx)$ fails to be in $L^2_{\mu}(\mathbb{R})$.
Hence we must have $\zeta>0$ and choose $d_1\in \mathbb{C}$ such that $\|f\|_{\mu,2}=1$.

When  $A(x)=-\frac{1}{\xi} (x-\langle x\rangle_{\mu})+\langle x\rangle_{D_\mu(f)}$, then following a similar reasoning as before, we can conclude that  $f(x)=d_2e^{-\frac{1}{2\zeta}x^2}e^{\frac{i}{2\xi}x^2}E_{\mu}(b^{\prime}x)$, where $b^{\prime}=\left(\frac{1}{\zeta}+\frac{i}{\xi}\right)\langle x\rangle_{f} +i \langle x\rangle_{D_{\mu}(f)}$, $\zeta > 0$ and $d_2\in \mathbb{C}$.

When
\begin{align*}
	A(x)=\left\{\begin{array}{cc}  \frac{1}{\xi} (x-\langle x\rangle_{f})+\langle x\rangle_{D_\mu(f)},& x\geq \langle x\rangle_{f},\qquad\qquad \\
		-\frac{1}{\xi} (x-\langle x\rangle_{\mu})+\langle x\rangle_{D_\mu(f)},& x <\langle x\rangle_{f},\\
	\end{array}\right.
\end{align*}
then, we must have
\begin{align*}
		f(x)=\left\{\begin{array}{cc} d_3 e^{-\frac{1}{2\zeta}x^2} e^{\frac{i}{2\xi}x^2} E_{\mu}(bx),& x\geq \langle x\rangle_{f}, \\
			d_4 e^{-\frac{1}{2\zeta}x^2} e^{-\frac{i}{2\xi}x^2} E_{\mu}(b^{\prime}x),& x <\langle x\rangle_{f}.\\
		\end{array}\right. 
\end{align*}
Finally, when
\begin{align*}
	A(x)=\left\{\begin{array}{cc}  -\frac{1}{\xi} (x-\langle x\rangle_{\mu})+\langle x\rangle_{D_\mu(f)},& x\geq \langle x\rangle_{f}, \\
	\frac{1}{\xi} (x-\langle x\rangle_{\mu})+\langle x\rangle_{D_\mu(f)},& x <\langle x\rangle_{f},
\end{array}\right.
\end{align*}
then
\begin{align*}
	f(x)=\left\{\begin{array}{cc} \quad d_5 e^{-\frac{1}{2\zeta}x^2} e^{-\frac{i}{2\xi}x^2} 	E_{\mu}(b^{\prime}x),& x\geq \langle x\rangle_{f}, \\
	d_6 e^{-\frac{1}{2\zeta}x^2} e^{ \frac{i}{2\xi}x^2} 	E_{\mu}(bx),& x <\langle x\rangle_{f}. \\
\end{array}\right. 
\end{align*}
Conversely, if $f$ has one of the forms (\ref{12}), (\ref{13}), (\ref{14}) or (\ref{15}), with $\|f\|_{\mu,2}=1$, then it is clear that (\ref{117}) and (\ref{118}) holds.
$\hfill\square$\\

Next, we proceed to prove the uncertainty inequality (\ref{16}) in Theorem \ref{MTHM1}.
\vspace{.2cm}

 \noindent{\bf{Proof of Theorem \ref{MTHM1}}}. 
	Assume that $\beta\notin\pi\mathbb{Z}$ and let $a=\langle x\rangle_{D^{\alpha}_{\mu}(f)}$ and $c=\langle x\rangle_{D^{\beta}_{\mu}(f)}$. Now using  Hausdorff-Young inequality (\ref{HYin}), we obtain
	\begin{align*}
	\Delta_{\mu,p}(D_{\mu}^{\beta}(f))&=\bigg(\int\limits_{-\infty}^{+\infty}\left|(x-c)D_{\mu}^{\beta}(f)(x)\right|^p \left|x\right|^{2\mu+1} dx\bigg)^{\frac{1}{p}}\\
		&=\bigg(\int\limits_{-\infty}^{+\infty}\left|(x-c)D_{\mu}^{\beta-\alpha}(D_{\mu}^{\alpha}f(u))(x)\right|^p \left|x\right|^{2\mu+1}dx\bigg)^{\frac{1}{p}}\qquad\qquad\\
		&=\frac{1}{|\sin(\beta-\alpha)|^{\mu+1}}\bigg(\int\limits_{-\infty}^{+\infty}\left|(x-c)D_{\mu}\left(e^{\frac{iu^2\cot(\beta-\alpha)}{2}}D_{\mu}^{\alpha}f(u)\right) \bigg(\frac{x}{\sin(\beta-\alpha)}\bigg)\right|^p |x|^{2\mu+1} dx\bigg)^{\frac{1}{p}}\\
		&=|\sin(\beta-\alpha)|^{-\mu+\frac{2(\mu+1)}{p}}\bigg(\int\limits_{-\infty}^{+\infty}\left|(x-c \csc(\beta-\alpha))D_{\mu}(g)(x)\right|^p |x|^{2\mu+1} dx\bigg)^{\frac{1}{p}}\\
		&=|\sin(\beta-\alpha)|^{-\mu+\frac{2(\mu+1)}{p}}\bigg(\int\limits_{-\infty}^{+\infty}\left|D_{\mu}((-iT_{\mu}-c^{\prime})g)(x)\right|^p |x|^{2\mu+1} dx\bigg)^{\frac{1}{p}}\\
		&\geq\frac{|\sin(\beta-\alpha)|^{-\mu+\frac{2(\mu+1)}{p}}}{\left(2^{\mu+1}\Gamma(\mu +1)\right)^{\frac{2}{p}-1}}\bigg(\int\limits_{-\infty}^{+\infty}\left|(-iT_{\mu}-c^{\prime})g(x)\right|^q |x|^{2\mu+1} dx\bigg)^{\frac{1}{q}},
	\end{align*}
	where $c^{\prime}=b\csc(\beta-\alpha),~~~g(x)=e^{\frac{ix^2\cot(\beta-\alpha)}{2}}D_{\mu}^{\alpha}f(x),~~x\in\mathbb{R}$. Thus by $\text{H\"{o}lder's}$ inequality,
	\begin{align*}
	&	\Delta_{\mu,p}(D_{\mu}^{\alpha}(f))\Delta_{\mu,p}(D_{\mu}^{\beta}(f))\\
		&\quad\geq \frac{|\sin(\beta-\alpha)|^{-\mu+\frac{2(\mu+1)}{p}}}{\left(2^{\mu+1}\Gamma(\mu +1)\right)^{\frac{2}{p}-1}}\bigg(\int\limits_{-\infty}^{+\infty}\left|(x-a)D_{\mu}^{\alpha}f(x)\right|^p |x|^{2\mu+1} dx\bigg)^{\frac{1}{p}} \bigg(\int\limits_{-\infty}^{+\infty}\left|(-iT_{\mu}-c^{\prime})g(x)\right|^q |x|^{2\mu+1} dx\bigg)^{\frac{1}{q}}\\
		&\quad\geq \frac{|\sin(\beta-\alpha)|^{-\mu+\frac{2(\mu+1)}{p}}}{\left(2^{\mu+1}\Gamma(\mu +1)\right)^{\frac{2}{p}-1}}\left|\int\limits_{-\infty}^{+\infty} e^{\frac{ix^2\cot(\beta-\alpha)}{2}}
		(x-a) D_{\mu}^{\alpha}f(x)\overline{(-iT_{\mu}-c^{\prime})g(x)}|x|^{2\mu+1}dx\right|.
	\end{align*}
Now using
	\begin{align*}
	T_{\mu}g(x)=ix\cot(\beta-\alpha)e^{\frac{ix^2 \cot(\beta-\alpha)}{2}}D_{\mu}^{\alpha}f(x)+e^{\frac{ix^2 \cot(\beta-\alpha)}{2}} T_{\mu}D_{\mu}^{\alpha}f(x),
	\end{align*}
we have
	\begin{align*}
	 \Delta_{\mu,p}(D_{\mu}^{\alpha}(f))\Delta_{\mu,p}(D_{\mu}^{\beta}(f))&\geq|\sin(\beta-\alpha)|^{-\mu+\frac{2(\mu+1)}{p}}\left| \text{{\bf{I}}}+\text{{\bf{II}}}\right|,
	\end{align*}
	where
	\begin{align}\label{k}
		\text{{\bf{I}}}=\int\limits_{-\infty}^{+\infty}x(x-a)~\cot(\beta-\alpha) \,|D_{\mu}^{\alpha}f(x)|^2 |x|^{2\mu+1}dx=\cot(\beta-\alpha) \, \Delta_{\mu,2}^2(D_{\mu}^{\alpha}f).
	\end{align}
and
\begin{align}
\text{{\bf{II}}}=\int\limits_{-\infty}^{+\infty} (x-a)D_{\mu}^{\alpha}f(x) \overline{(-iT_{\mu}-c^{\prime})D_{\mu}^{\alpha}f(x)}|x|^{2\mu+1}dx.
\end{align}
	{\it{Case 1:}} Assume that $\alpha\neq n\pi$, for some $n\in\mathbb{Z}$, then
	\begin{align}
		D_{\mu}^{\alpha}f(x)=\frac{1}{|\sin\alpha|^{\mu+1}} e^{\frac{ix^2 \cot\alpha}{2}}D_{\mu}(h)\left(\frac{x}{\sin\alpha}\right), \, \text{where} \, h(x)=e^{\frac{ix^2\cot\alpha}{2}}f(x).
	\end{align}
Now,
	\begin{align*}
		T_{\mu}(D_{\mu}f(x))=\frac{1}{|\sin\alpha|^{\mu+1}}e^{i\frac{x^2\cot\alpha}{2}}\bigg[ix \cot\alpha \, D_{\mu}(h)\left(\frac{x}{\sin\alpha}\right)+T_{\mu}(D_{\mu}h)\left(\frac{x}{\sin\alpha}\right)\bigg].
	\end{align*}
	Therefore, 
	\begin{align*}
		\text{{\bf{II}}}&=\frac{\cot\alpha}{|\sin\alpha|^{(\mu+1)2}}\int\limits_{-\infty}^{+\infty} x(x-a)\left|D_{\mu}(h)\left(\frac{x}{\sin\alpha}\right)\right|^2 |x|^{2\mu+1}dx\qquad\qquad\qquad\\
		&\qquad+\frac{1}{|\sin\alpha|^{2(\mu+1)}}\int\limits_{-\infty}^{+\infty} (x-a)D_{\mu}(h)\left(\frac{x}{\sin\alpha}\right)(-iT_{\mu}-c^{\prime})\left[D_{\mu}(h)\left(\frac{x}{\sin\alpha}\right)\right](x)~~|x|^{2\mu+1} dx\\
		&=\text{{\bf{II}}}_1+\text{{\bf{II}}}_2.
	\end{align*}
	Now, 
	\begin{align*}
		\text{{\bf{II}}}_1&=\sin\alpha \, \cos\alpha\, \int\limits_{-\infty}^{+\infty} x\left(x-\frac{a}{\sin\alpha}\right)\left|D_{\mu}(h)(x)\right|^2 |x|^{2\mu+1} dx\\
		&=\sin\alpha\, \cos\alpha \, \int\limits_{-\infty}^{+\infty} x^2\left|D_{\mu}(h)(x)\right|^2 \left|x\right|^{2\mu+1}dx-a \cos\alpha\,  \int\limits_{-\infty}^{+\infty} x\left|D_{\mu}(h)(x)\right|^2 \left|x\right|^{2\mu+1}dx\\
		&=\sin\alpha\, \cos\alpha\, \int\limits_{-\infty}^{+\infty} (x-\langle x\rangle_{D_{\mu}(h)})^2\left|D_{\mu}(h)(x)\right|^2 \left|x\right|^{2\mu+1}dx+\langle x\rangle_{D_{\mu}(h)}\,\cos\alpha \,\left(\sin\alpha \,  \langle x\rangle_{D_{\mu}(h)}-a\right).
	\end{align*}
	By Lemma (\ref{lem2}) and Lemma (\ref{lem3}), we have  $\sin\alpha \,\langle x\rangle_{D_{\mu}(h)}=a$. So
	\begin{align*}
		\text{{\bf{II}}}_1=\sin\alpha \, \cos\alpha \Delta_{\mu,2}^2(D_{\mu}(h)).
	\end{align*}
Also,
	\begin{align*}
	\text{{\bf{II}}}_2&=\frac{1}{\left|\sin\alpha\right|^{(\mu+1)2}}\times\frac{1}{\sin\alpha}\int\limits_{-\infty}^{+\infty} (x-a) D_{\mu}(h)\left(\frac{x}{\sin\alpha}\right)\left(-iT_{\mu}-\sin\alpha \, c^{\prime}\right)D_{\mu}(h)\left(\frac{x}{\sin\alpha}\right) |x|^{2\mu+1}dx\qquad\qquad\\
		&=\int\limits_{-\infty}^{+\infty} \left(x-\frac{a}{\sin\alpha}\right)D_{\mu}(h)(x)\overline{(-iT_{\mu}-\sin\alpha\, c^{\prime})D_{\mu}(h)(x)} |x|^{2\mu+1}dx\\
		&=\int\limits_{-\infty}^{+\infty} xD_{\mu}(h)(x)\overline{(-iT_{\mu})D_{\mu}(h)(x)} |x|^{2\mu+1}dx+ac^{\prime}\int\limits_{-\infty}^{+\infty}\left|D_{\mu}(h)(x)\right|^2 |x|^{2\mu+1} dx\\
		&\qquad-\sin\alpha \, c^{\prime}\int\limits_{-\infty}^{+\infty} x\left|D_{\mu}(h)(x)\right|^2 |x|^{2\mu+1} dx-\frac{a}{\sin\alpha}\int\limits_{-\infty}^{+\infty}D_{\mu}(h)(x)\overline{(-iT_{\mu}(D_{\mu}(h)(x)))}\left|x\right|^{2\mu+1} dx.
	\end{align*}
Clearly $\int\limits_{-\infty}^{+\infty}\left|D_{\mu}(h)(x)\right|^2 |x|^{2\mu+1} dx = 1$, and from Lemma \ref{lem2} and Lemma \ref{lem3}, we have $a = \sin \alpha \langle x \rangle_{D_\mu(h)}$ and
\begin{align*}
		\int\limits_{-\infty}^{+\infty}D_{\mu}(h)(x)\overline{(-iT_{\mu}(D_{\mu}(h)(x)))}\left|x\right|^{2\mu+1} dx
		&=\int\limits_{-\infty}^{+\infty} h(-x)\overline{xh(-x)}\left|x\right|^{2\mu+1}dx\\
		&=-\int\limits_{-\infty}^{+\infty} x\left|f(x)\right|^2 |x|^{2\mu+1} dx=-\langle x\rangle_f.
	\end{align*}
  	Substituting this back into  $\text{{\bf{II}}}_2$ and  using (\ref{re4}) we have
	\begin{align*}
		\text{{\bf{II}}}_2&=\int\limits_{-\infty}^{+\infty} iT_{\mu}h(x)\overline{xh(x)}\left|x\right|^{2\mu+1}dx + \frac{a}{\sin\alpha}\langle x\rangle_f\\
		&=-\int\limits_{-\infty}^{+\infty} \cot\alpha \, x^2 \left|f(x)\right|^2\left|x\right|^{2\mu+1}dx + i\int\limits_{-\infty}^{+\infty} T_{\mu}f(x){\overline{xf(x)}}\left|x\right|^{2\mu+1}dx+\frac{a}{\sin\alpha}\langle x\rangle_f\\
		&=-\cot\alpha \, \int\limits_{-\infty}^{+\infty}(x-\langle x \rangle_f )^2\left|f(x)\right|^2\left|x\right|^{2\mu+1}dx-\cot\alpha\langle x\rangle _f^2\\
		&\quad\quad\quad\quad\quad\quad\quad\quad + i\left[- \mathcal{A} + i \left(Cov_{\mu}(f)+\langle x \rangle_f \langle x \rangle_{D_{\mu(f)}}\right) \right]+\frac{a}{\sin\alpha}\langle x \rangle_f\\
		&=-\cot\alpha \Delta_{\mu,2}^2 (f)-Cov_{\mu}(f)-i\mathcal{A}+\frac{\langle x\rangle_{f}}{2}\left[-\cos\alpha \langle x\rangle _f+a \right]- \langle x\rangle _f \langle x\rangle_{D_{\mu}(f)}\\
		&=-\cot\alpha \, \Delta_{\mu,2}^2(f)-Cov_{\mu}(f)-i\mathcal{A},
	\end{align*}
where $\mathcal{A}$ is given in (\ref{uyt}). Therefore,
	\begin{align*}
		\text{{\bf{II}}}=\text{{\bf{II}}}_1+\text{{\bf{II}}}_2=\left(\frac{\cos^3\alpha-\cos\alpha}{\sin\alpha}\right)\Delta_{\mu,2}^2(f)+(2\cos^2\alpha-1)Cov_{\mu}(f)+\sin\alpha \,\cos\alpha \Delta_{\mu,2}^2(D_{\mu}(f))-i\mathcal{A}.
	\end{align*}
	So, we have
	\begin{align*}
		\text{{\bf{I}}}+\text{{\bf{II}}}&=\left(\cot(\beta-\alpha)\cos^2\alpha+\frac{\cos^3\alpha-\cos\alpha}{\sin\alpha}\right)\Delta_{\mu,2}^2(f)+\bigg(2\sin\alpha\, \cos\alpha\cot(\beta-\alpha) + 2\cos^2\alpha-1\bigg)\\
		&\qquad\times Cov_{\mu}(f)+\bigg(\sin^2\alpha \cot(\beta-\alpha) + \sin\alpha\, \cos\alpha\bigg)\Delta_{\mu,2}^2(D_{\mu}(f))-i\mathcal{A}\\
		&=\left(\frac{\cos\alpha\, \cos\beta}{\sin(\beta-\alpha)}\right)\Delta_{\mu,2}^2(f)+\left(\frac{\sin(\alpha+\beta)}{\sin(\beta-\alpha)}\right)Cov_{\mu}(f)+\left(\frac{\sin\alpha\, \sin\beta}{\sin(\beta-\alpha)}\right)\Delta_{\mu,2}^2(D_{\mu}(f))-i\mathcal{A}.\\
	\end{align*}
	Thus, \begin{align*}
		\left|\text{{\bf{I}}}+\text{{\bf{II}}}\right|&=\left\{\left[\left(\frac{\cos\alpha\cos\beta}{\sin(\beta-\alpha)}\right)\Delta_{\mu,2}^2(f)+\left(\frac{\sin(\alpha+\beta)}{\sin(\beta-\alpha)}\right)Cov_{\mu}(f)+\left(\frac{\sin\alpha\sin\beta}{\sin(\beta-\alpha)}\right)\Delta_{\mu,2}^2(D_{\mu}(f))\right]^2 +\mathcal{A}^2 \right\}^{\frac{1}{2}}.\qquad\qquad
	\end{align*}
	Now,
	\begin{align*}
		\Delta_{\mu, p}^2( D_{\mu}^{\alpha}f)\Delta_{\mu,p}^2(D_{\mu}^{\beta}f)&\geq\frac{|\sin(\beta-\alpha)|^{-2\mu+\frac{4(\mu+1)}{p}}}{\left(2^{\mu+1}\Gamma(\mu +1)\right)^{2(\frac{2}{p}-1)}}\left|\text{I}+\text{II}\right|^2\\
		&=\frac{|\sin(\beta-\alpha)|^{-2\mu+\frac{4(\mu+1)}{p}}}{\left(2^{\mu+1}\Gamma(\mu +1)\right)^{2(\frac{2}{p}-1)}}\biggl\{\left[\left(\frac{\cos\alpha\cos\beta}{\sin(\beta-\alpha)}\right)\Delta_{\mu,2}^2(f)+\left(\frac{\sin(\alpha+\beta)}{\sin(\beta-\alpha)}\right)Cov_{\mu}(f)\right.\\
		&\qquad\qquad\left.+\left(\frac{\sin\alpha\sin\beta}{\sin(\beta-\alpha)}\right)\Delta_{\mu,2}^2(D_{\mu}(f))\right]^2+\mathcal{A}^2\biggr\}.
	\end{align*}
	Consequently, we obtain (\ref{16}).\\
	{\it{Case 2:~}}~ Assume that $\alpha= 2n\pi$, for some $n\in\mathbb{Z}$. Then from (\ref{k}), we have $\text{{\bf{I}}}=\cot(\beta)\,\Delta_{\mu,2}^2(f)$ and 
	\begin{align*}
		\text{{\bf{II}}}&=\int\limits_{-\infty}^{+\infty} (x-a)f(x)\overline{\left(-iT_{\mu}- b^{\prime}\right)f(x)} |x|^{2\mu+1}dx\\
		&=\int\limits_{-\infty}^{+\infty} xf(x) \overline{\left(-iT_{\mu}\right)f(x)}|x|^{2\mu+1} dx+ab^{\prime}\int\limits_{-\infty}^{+\infty} |f(x)|^2 |x|^{2\mu+1}dx\\
		&-a\int\limits_{-\infty}^{+\infty} f(x)\overline{\left(-iT_{\mu}\right)f(x)} |x|^{2\mu+1}dx-b^{\prime}\int\limits_{-\infty}^{+\infty} x|f(x)|^2 |x|^{2\mu+1}dx\\
		&=\int\limits_{-\infty}^{+\infty}xf(x)\overline{\left(-iT_{\mu}\right)f(x)}|x|^{2\mu+1} dx-\langle x\rangle_f \langle x\rangle_{D_{\mu}(f)}\\
		&=Cov_{\mu}(f)+i\mathcal{A}.
	\end{align*}
Therefore,
	\begin{align*}
		\Delta_{\mu, p}^2( D_{\mu}^{\alpha}f)\Delta_{\mu,p}^2(D_{\mu}^{\beta}f)&\geq\frac{|\sin(\beta-\alpha)|^{-2\mu+\frac{4(\mu+1)}{p}}}{\left(2^{\mu+1}\Gamma(\mu +1)\right)^{2(\frac{2}{p}-1)}}\left(\mathcal{A}^2+\left(\cot\beta \, \Delta_{\mu,2}^2+Cov_{\mu}(f)\right)^2\right),
	\end{align*}
thus (\ref{16}) follows.\\
{\it{Case 3:~}}~ Assume that $\alpha= (2n+1)\pi$, for some $n\in\mathbb{Z}$. Arguing as in the previous case, we obtain 
	\begin{align*}
	\Delta_{\mu, p}^2( D_{\mu}^{\alpha}f)\Delta_{\mu,p}^2(D_{\mu}^{\beta}f)&\geq\frac{|\sin(\beta-\alpha)|^{-2\mu+\frac{4(\mu+1)}{p}}}{\left(2^{\mu+1}\Gamma(\mu +1)\right)^{2(\frac{2}{p}-1)}}\left(\mathcal{A}^2+\left(-\cot\beta \, \Delta_{\mu,2}^2+Cov_{\mu}(f)\right)^2\right).
\end{align*}
As a result, we get (\ref{16}).
$\hfill\square$\\

Now we prove Theorem \ref{MTHM2}.
\vspace{.2cm}

\noindent{\bf{Proof of Theorem \ref{MTHM2}}}. 
By using Lemma \ref{lem2} and Theorem \ref{fe}, we have
\begin{align} \label{s123}
	\notag	\Delta_{\mu, 2}^2( D_{\mu}^{\alpha}f)\Delta_{\mu,2}^2(D_{\mu}^{\beta}f)
	&=\sin^2(\alpha-\beta) \,\left(\Delta_{\mu,2}^2(f)\Delta_{\mu,2}^2(D_{\mu}(f))-Cov_{\mu}^2(f)\right)\\ \notag
	&+\biggl(\cos\alpha\, \cos\beta\, \Delta_{\mu, 2}^2(f)+\sin\alpha\, \sin\beta\, \Delta_{\mu, 2}^2(D_{\mu}(f))+\sin(\alpha+\beta)\, Cov_{\mu}(f)\biggr)^2\\ \notag
	\geq&\sin^2(\alpha-\beta)\left(\left\{\left({\mu}+\frac{1}{2}\right)\left(\|f_e\|^2_{\mu,2}-\|f_0\|^2_{\mu,2}\right)+\frac{1}{2}\right\}^2+COV_{\mu}^2(f)-Cov_{\mu}^2(f)\right)\\
	&+\biggl(\cos\alpha\, \cos\beta\, \Delta_{\mu, 2}^2(f)+\sin\alpha\, \sin\beta\, \Delta_{\mu, 2}^2(D_{\mu}(f))+\sin(\alpha+\beta)\, Cov_{\mu}(f)\biggr)^2.
\end{align}
Notice that equality holds in (\ref{s123}) if and only if the equality holds in (\ref{110}).

$\hfill\square$\\

Now, it remains to prove the conditions under which equality holds in Theorem \ref{MTHM1}.
\vspace{.2cm}

\noindent{\bf{Proof of the case when the lower bound is attained in Theorem \ref{MTHM1}}}. Since we are applying the Hausdorff-Young inequality to prove equation (\ref{16}), it follows that the condition $p=2$ is necessary for equality to hold in (\ref{16}). If $p=2$, then equality holds in (\ref{16}) if and only if equality holds in (\ref{17}) with $COV_\mu(f) = Cov_\mu(f)$. Moreover, $COV_\mu(f) = Cov_\mu(f)$ if and only if $(x-\langle x \rangle_{f})$ and  $T_{\mu}(\varphi)(x) + \frac{\mu+\frac{1}{2}}{x}\bigg(\frac{\sin[\varphi(x)-\varphi(-x)] \rho(-x)}{\rho(x)}-(\varphi(x)-\varphi(-x))\bigg)
-\langle x\rangle_{D_\mu(f)}$  have the same or opposite signs (see Subsection \ref{sub2.3}).

Thus, equality holds in (\ref{16}) if and only if $f(x)=\rho(x)e^{i\varphi(x)}$ satisfies either equations (\ref{117}) and (\ref{38}), or equations (\ref{117}) and (\ref{39}). Therefore, $f$ must be of the form (\ref{12}) or (\ref{13}). This completes the proof of Theorem \ref{MTHM1}.
$\hfill\square$

\vspace{.3cm}
\begin{remark}
	(1) The constants $d_1, d_2$ in Theorem \ref{THM1} can be chosen such that (see Lemma 4.3 in \cite{Rosler})  $d_1 = e^{i\theta} \left(\zeta^{\mu +1} \Gamma(\mu +1)  e^{- Re\left(\frac{b^2\zeta}{2}\right)}E_\mu(\frac{|b|^2\zeta}{2})\right)^{-\frac{1}{2}}$and $ d_2 = e^{i\theta} \left(\zeta^{\mu +1} \Gamma(\mu +1)  e^{- Re\left(\frac{b'^2\zeta}{2}\right)}E_\mu(\frac{|b'|^2\zeta}{2})\right)^{-\frac{1}{2}}$.\\
	(2) For $\mu > -\frac{1}{2}$, the weighted measure $|x|^{2\mu +1} dx$ is not translation-invariant. As a result, Theorem \ref{MTHM1} cannot be recovered by first proving it for the cases where $\langle x \rangle_{f} = 0$,$\langle x \rangle_{D_\mu(f)} = 0$ as in Chapter 6 of \cite{1} or \cite{Chen}.\\
	(3) If the function $f$ is of the form given in (\ref{14}) or (\ref{15}), it is clear from the proof of Theorem \ref{THM1} that $(x-\langle x \rangle_{f})$ and  $T_{\mu}(\varphi)(x) + \frac{\mu+\frac{1}{2}}{x}\bigg(\frac{\sin[\varphi(x)-\varphi(-x)] \rho(-x)}{\rho(x)}-(\varphi(x)-\varphi(-x))\bigg)
	-\langle x\rangle_{D_\mu(f)}$ do not have the same or opposite signs. Hence, it follows that $COV_\mu(f) > cov_\mu(f)$ (see Subsection \ref{sub2.3}).
\end{remark}

\section*{Acknowledgments}
\baselineskip=11pt

The first and the second author wishes to thank the Ministry of Human Resource Development, India for the  research fellowship and Indian Institute of Technology Guwahati, India for the support provided during the period of this work.

\end{document}